\documentclass[12pt]{article}
\usepackage{amsmath,fullpage,amssymb,hyperref,amsthm}
\usepackage{tikz}

\newtheorem{theorem}{Theorem} 

\newtheorem{corollary}[theorem]{Corollary}
   
\newcommand\card[1]{\left|#1\right|}
\newcommand\si{\sigma}

\def\S{\mathfrak{S}}
\DeclareMathOperator\Av{Av}
\DeclareMathOperator\wt{wt}
\newcommand{\ul}{\underline}

\title{Consecutive patterns in circular permutations}
\author{Sergi Elizalde\\
Department of Mathematics\\
Dartmouth College\\
Hanover, NH 03755-3551, USA\\
\texttt{sergi.elizalde@dartmouth.edu}
\and Bruce Sagan\\
Department of Mathematics\\
Michigan State University\\
East Lansing, MI 48824-1027, USA\\
\texttt{bsagan@msu.edu}
}
\date{}

\begin{document}

\maketitle

\begin{abstract} 
In their study of cyclic pattern containment, Domagalski et al.~\cite{DLMSSS} conjecture differential equations for the generating functions of circular permutations avoiding consecutive patterns of length~3. In this note, we prove and significantly generalize these conjectures. 
We show that, for every consecutive pattern $\si$ beginning with $1$, the bivariate generating function counting occurrences of $\si$ in circular permutations can be obtained from the generating function counting occurrences of $\si$ in (linear) permutations. This includes all the patterns for which the latter generating function is known.
\end{abstract}

\section{Introduction}

\subsection{Pattern containment and avoidance}

Given two permutations $\pi=\pi_1\pi_2\dots\pi_n\in\S_n$ and  $\si_1\si_2\dots\si_m\in\S_m$, 
we say that $\pi$ {\em contains} the {\em classical pattern} $\si_1\si_2\dots\si_m$
if there is a subsequence $\pi_{i_1}\pi_{i_2}\dots\pi_{i_m}$, where $1\le i_1<i_2<\dots<i_m\le n$, that is {\em order-isomorphic} to $\si_1\si_2\dots\si_m$, meaning that $\pi_{i_a}<\pi_{i_b}$ if and only if $\si_a<\si_b$ for all $1\le a,b\le m$.
We say that $\pi$ {\em contains} the {\em consecutive pattern} $\ul{\si_1\si_2\dots\si_m}$
(we underline consecutive patterns to distinguish them from classical ones)
if there is a consecutive subsequence $\pi_{i}\pi_{i+1}\dots\pi_{i+m-1}$, where $1\le i\le n-m+1$, that is order-isomorphic to $\si_1\si_2\dots\si_m$. Writing $\si=\ul{\si_1\si_2\dots\si_m}$, such a subsequence is called an {\em occurrence} of $\si$ in $\pi$, and the number of such occurrences is denoted by $o_\si(\pi)$.
If $\pi$ does not contain a pattern, we say that $\pi$ {\em avoids} the pattern.

While classical patterns have been extensively studied for over half a century~\cite{SS}, the systematic study of consecutive patterns in permutations started only two decades ago~\cite{EliNoy1}. 
Around the same time, {\em vincular patterns}, which generalize classical and consecutive patterns by requiring certain positions to be adjacent, were introduced in~\cite{BS} by the name {\em generalized patterns}.

Several authors~\cite{Callan,Vella,DLMSSS} have considered another variation of the notion of pattern containment, where one allows the patterns to wrap around from the end of the permutation to the beginning. This variation has received different names in the literature.

Denote the set of rotations of a permutation $\pi=\pi_1\pi_2\dots\pi_n\in\S_n$ by
$$[\pi]=\{\pi_1\pi_2\dots\pi_n,\pi_2\dots\pi_n\pi_2,\dots,\pi_n\pi_1\dots\pi_{n-1}\}.$$
(These are sometimes called {\em horizontal rotations}~\cite{EliRoi1,EliRoi2} to distinguish them from {\em vertical rotations}, which are obtained when shifting the values of the entries instead of their positions.)
This set is called a {\em circular permutation} by Callan~\cite{Callan}, a {\em cyclic arrangement} by Vella~\cite{Vella}, and a {\em cyclic permutation} by Domagalski et al.~\cite{DLMSSS}. In this note, we will use the term {\em circular permutations} for such sets, in order to distinguish them from permutations whose cycle structure consists on a single cycle, which, incidentally, are also quite interesting from a pattern avoidance perspective~\cite{ArcEli,Elicont}.
As in~\cite{DLMSSS}, we denote by $[\S_n]$ the set of circular permutations of length $n$, namely, the set of equivalence classes of permutations in $\S_n$ under rotation. Permutations in $\S_n$ may be called {\em linear} permutations when we want to distinguish them from circular ones.

For any pattern $\si$, whether classical ($\si=\si_1\si_2\dots\si_m$) or consecutive ($\si=\ul{\si_1\si_2\dots\si_m}$), we say that the circular permutation $[\pi]$ {\em contains} $\si$ if there is some rotation $\pi'\in[\pi]$ that contains $\si$ according to the above definition for linear permutations; otherwise, we say that $[\pi]$ {\em avoids} $\si$. Following~\cite{DLMSSS}, we denote by $\Av_n[\si]$ the set of circular permutations in $[\S_n]$ that avoid $\si$.

While the enumeration of circular permutations avoiding classical patterns of length 3 is trivial, Vella~\cite{Vella} and Callan~\cite{Callan} determined $\card{\Av_n[\si]}$ for all classical patterns $\si$ of length $4$. In a recent preprint, Domagalski et al.~\cite{DLMSSS} enumerate circular permutations avoiding any subset of patterns of length~$4$. All the formulas that have been obtained so far have simple expressions involving binomial coefficients, linear terms, powers of 2, constants, and Fibonacci numbers.

On the other hand, the enumeration of circular permutations avoiding consecutive patterns has not yet been explored. This is left as an open problem in~\cite{DLMSSS}, where some conjectures are made \cite[Conjecture 6.4]{DLMSSS} in the special case of patterns of length~3, in the form of differential equations hypothetically satisfied by the generating functions
for the permutations avoiding them. These are stated as Equations~\eqref{eq:conj1} and~\eqref{eq:conj2} below.
The goal of this note is to prove these conjectures, and to generalize them in two directions. On the one hand, we enumerate not only circular permutations avoiding each pattern, but also circular permutations with any given number of occurrences of the pattern. On the other hand, we extend the results to other consecutive patterns, namely all of those for which the generating function tracking their occurrences in linear permutations is currently known.

In the rest of the paper, we let $\si=\ul{\si_1\si_2\dots\si_m}$ be a consecutive pattern.
We define an {\em occurrence} of the consecutive pattern $\si$ in a circular permutation $[\pi]\in[\S_n]$ to be a subsequence $\pi_{i}\pi_{i+1}\dots\pi_{i+k-1}$
or $\pi_i\pi_{i+1}\dots\pi_n\pi_1\pi_2\dots\pi_{k-n+i-1}$ (i.e., allowed to wrap around), where $1\le i\le n$, that is order-isomorphic to $\si_1\si_2\dots\si_m$. We denote by $c_\si[\pi]$ the number of occurrences of $\si$ in $[\pi]$. This number is well defined, in the sense that it does not depend on the chosen representative of $[\pi]$, since rotating $\pi$ simply changes the positions of the occurrences of $\si$, but not the actual subsequences or how many there are.
Note also that $c_\si[\pi]=0$ precisely if $[\pi]$ avoids $\si$. For example, $c_{\ul{132}}[25314]=2$, since $253$ and $142$ are occurrences of $\si$ in $[\pi]$. On the other hand, $c_{\ul{132}}[24531]=0$, so $[24531]\in\Av_5[\ul{132}]$.

\subsection{Generating functions}

We denote by
$$P_\si(u,z)=\sum_{n\ge0}\sum_{\pi\in\S_n} u^{o_\si(\pi)}\frac{z^n}{n!}$$
the exponential generating function counting occurrences of a consecutive pattern $\si$ in linear permutations, and let $\omega_\si(u,z)=1/P_\si(u,z)$. 
Formulas and differential equations for $P_\si(u,z)$ and $\omega_\si(u,z)$, for various patterns $\si$, have been given in~\cite{EliNoy1,EliNoy2}, see also~\cite{Kit,MR,LR,DK} for related work.

In this paper, we are interested in the analogues for circular permutations of these generating functions.
Let
\begin{equation}\label{def:C}
C_{\si}(u,z)=\sum_{n\ge0} \sum_{[\pi]\in[\S_n]} u^{c_\si[\pi]}\frac{z^n}{n!}
\end{equation}
be the exponential generating function counting occurrences of $\si$ in circular permutations, and note that
$$C_{\si}(0,z)=\sum_{n\ge0} \card{\Av_n[\si]}\frac{z^n}{n!}.$$

As in the case of consecutive patterns in linear permutations, letting $\si^r=\ul{\si_m\dots\si_2\si_1}$ and $\si^c=\ul{(m+1-\si_1)(m+1-\si_2)\dots(m+1-\si_m)}$, it is clear that $$C_{\si}(u,z)=C_{\si^c}(u,z)=C_{\si^r}(u,z)=C_{\si^{rc}}(u,z),$$
since occurrences of $\si$ in $[\pi]$ correspond to occurrences of $\si^r$ in $[\pi^r]$, and to occurrences of $\si^c$ in $[\pi^c]$. 
For example, for patterns of length~3, we have $C_{\ul{123}}(u,z)=C_{\ul{321}}(u,z)$
and  $C_{\ul{132}}(u,z)=C_{\ul{312}}(u,z)=C_{\ul{213}}(u,z)=C_{\ul{231}}(u,z)$.

For a function $F(u,z)$, we will use $F'(u,z)$ to denote its partial derivative with respect to the variable~$z$.

\section{Counting patterns in circular permutations}

Our central result relates consecutive patterns in the circular case with those in the linear case. 
The requirement $\sigma_1=1$ can be replaced, by the above symmetries, with any of $\sigma_1=m$, $\sigma_m=1$, or $\sigma_m=m$, where $m$ is the length of $\sigma$.

For the purposes of the proof we will let a {\em permutation} be any linear or cyclic ordering of a finite set of positive integers.
Any set of circular permutations  
$[\Pi]=\{[\pi^{(1)}],[\pi^{(2)}],\ldots,[\pi^{(k)}]\}$
will be given {\em weight}
$$
\wt [\Pi]= u^{c_{\si}[\pi^{(1)}]} \cdot 
u^{c_{\si}[\pi^{(2)}]} \cdots u^{c_{\si}[\pi^{(k)}]},
$$
and any linear permutation $\pi$  will be given weight
$\wt \pi=u^{o_{\si}}(\pi)$.  Finally, the {\em left-right minima} of $\pi=\pi_1\pi_2\ldots\pi_n$ are the elements $\pi_i$ such that
$$
\pi_i = \min\{\pi_1,\pi_2,\ldots,\pi_i\}.
$$
These elements give rise to the {\em left-right minima factorization} of $\pi$ which is 
\begin{equation}
\label{lrf}
\pi=\pi^{(1)}\pi^{(2)}\ldots\pi^{(k)}
\end{equation}
where $\pi^{(i)}$ is the factor (consecutive subword) of $\pi$ starting at the $i$th left-right minimum and ending just before the $(i+1)$st.

\begin{theorem}\label{thm:main}
Let $\si=\ul{\si_1\si_2\dots\si_m}$ be a consecutive pattern with $\si_1=1$. Then
$$C_\si(u,z)=1+\ln P_\si(u,z).$$
\end{theorem}

\begin{proof}
Exponentiating the equation in the statement of the theorem, it suffices to prove that
$$
P_\si(u,z) = e^{C_\si(u,z)-1}.
$$
By the Exponential Formula (see Theorem 4.5.1 in Sagan's book~\cite{sagan}), it suffices to show that there is a bijection $\phi$ between permutations $\pi\in\S_n$ and sets of circular permutations
$[\Pi]= \{[\pi^{(1)}],[\pi^{(2)}],\ldots,[\pi^{(k)}]\}$
such that
\begin{enumerate}
    \item[(a)] $\biguplus_{i=1}^{k} \pi^{(i)}=\{1,2,\ldots,n\}$, the union being of the underlying sets of the $\pi^{(i)}$, and
    \item[(b)] $\wt\pi=\wt[\Pi]$.
\end{enumerate}

Define
$$
\phi(\pi) = \{[\pi^{(1)}],[\pi^{(2)}],\ldots,[\pi^{(k)}]\}
$$
where the $\pi^{(i)}$ are the factors in~\eqref{lrf}.  Then (a) holds because every element of $\{1,2,\dots,n\}$ must appear in exactly one of the factors of the factorization.  To prove (b), let us show that any occurrence of $\si$ in $\pi$ is entirely contained in one of the $\pi^{(i)}$.  Indeed, if the occurrence overlaps two or more factors, then the left-right minimum of the second factor is smaller than the first element of the occurrence.  This contradicts the fact that $\si$ begins with~$1$.

To show $\phi$ is bijective, we construct its inverse.
Given $[\Pi]$, rotate each circular permutation so that 
$\pi^{(i)}$ starts with its minimum element.  Then concatenate these linear permutations in order of decreasing first element to form $\pi$.  It is easy to check that this describes the inverse of~$\phi$.
\end{proof}

\section{Applications to specific patterns}

Expressions for $P_\sigma=P_\sigma(u,z)$ are known for certain consecutive patterns $\sigma$, often in the form of differential equations satisfied by its reciprocal $\omega_\si=1/P_\si$. In fact, up to symmetry, all the patterns $\sigma$ for which explicit differential equations have been found so far satisfy $\sigma_1=1$. Thus, Theorem~\ref{thm:main} can be applied to these patterns to deduce an expression for $C_\si=C_\si(u,z)$. 

Restating Theorem~\ref{thm:main} to relate $C_\si$ and $\omega_\si$, we have $C_\si=1-\ln \omega_\si$, from where $C'_\si=-\omega'_\si/\omega_\si$, and 
\begin{equation}\label{eq:omC} \omega_\si=e^{1-C_\si}. \end{equation}
In some cases, this relation allows us to obtain differential equations directly in terms of $C_\si$, as we will see below.

\subsection{Monotone patterns}

It is proved in \cite[Theorem 3.1]{EliNoy1} (see also \cite[Theorem 2.1]{EliNoy2}) that, for $\si=\ul{12\dots m}$ with $m\ge3$, the function $\omega_\si=\omega_\si(u,z)$ satisfies the differential equation
\begin{equation}\label{eq:mon}
\omega_\si^{(m-1)}+(1-u)(\omega_\si^{(m-2)}+\dots+\omega_\si'+\omega_\si)=0
\end{equation}
with initial conditions $\omega_\si(u,0)=1$, $\omega_\si'(u,0)=-1$, and $\omega_\si^{(i)}(u,0)=0$ for $2\le i\le m-2$. 
In \cite[Theorem 2.4]{EliNoy2}, similar differential equations are given for $\omega_\sigma$ whenever $\sigma$ is a so-called {\em chain pattern} (see  \cite[Definition 2.2]{EliNoy2}). Chain patterns generalize monotone patterns, but they still satisfy $\sigma_1=1$ (up to symmetry), as shown in \cite[Lemma 2.3]{EliNoy2}. Thus, for all such patterns $\sigma$, Theorem~\ref{thm:main} can be used to determine $C_\si=1-\ln \omega_\si$.

It is possible to rewrite~\eqref{eq:mon} as a differential equation for $C_\sigma$ using the identity~\eqref{eq:omC}. For example, when $m=3$, we obtain the following.

\begin{corollary}
Let $D=D_{\ul{123}}(u,z)=C_{\ul{123}}'(u,z)$. Then $D$ satisfies the differential equation
\begin{equation}\label{eq:D123} D'=D^2+(u-1)(D-1) 
\end{equation}
with initial condition $D(u,0)=1$.
An explicit expression is given by
$$D_{\ul{123}}(u,z)=\frac{1}{2}\left(1-u-
\tanh \left( \frac{z\sqrt {{u}^{2}+2u-3}}{2}-{\rm arctanh} \left({
\frac {u+1}{\sqrt {{u}^{2}+2\,u-3}}}\right) \right) \sqrt {{u}^{2}+2\,
u-3}\right)
,$$
which, for $u=0$, simplifies to 
$$D_{\ul{123}}(0,z)=\frac{1}{2}+\frac{\sqrt{3}}{2}\tan\left(\frac{\sqrt{3}}{2}z+\frac{\pi}{6}\right).$$
\end{corollary}

\begin{proof}
Differentiating Equation~\eqref{eq:omC}, we get 
$\omega'_\si=-C'_\si \, e^{1-C_\si}$ and $\omega''_\si=\left(-C''_\si+(C'_\si)^2\right) e^{1-C_\si}$. Substituting these expressions into Equation~\eqref{eq:mon} for $m=3$, and dividing both sides by $e^{1-C_\si}$, we obtain Equation~\eqref{eq:D123}.
\end{proof}

Setting $u=0$ in Equation~\eqref{eq:D123} gives 
\begin{equation}\label{eq:conj1} D'_{\ul{123}}(0,z)=D_{\ul{123}}(0,z)^2-D_{\ul{123}}(0,z)+1, 
\end{equation}
proving part 1 of \cite[Conjecture 6.4]{DLMSSS}\footnote{Precisely speaking, the statement in part 1 of \cite[Conjecture 6.4]{DLMSSS} is slightly inaccurate, since
the equation that it gives is the one satisfied by $C_{\ul{123}}'(0,z)$, rather than by $C_{\ul{123}}(0,z)$.}. 
For $m=4$, a similar computation yields the following.

\begin{corollary}
Let $D=D_{\ul{1234}}(u,z)=C_{\ul{1234}}'(u,z)$. Then $D$ satisfies the differential equation
\begin{equation}\label{eq:D1234} D''=3D'D-D^3+(u-1)(D'-D^2+D-1) 
\end{equation}
with initial conditions $D(u,0)=1$, $D'(u,0)=1$.
For $u=0$, an explicit expression is given by
$$D_{\ul{1234}}(0,z)=\frac{\cos z+\sin z+e^{-z}}{\cos z-\sin z+e^{-z}}.$$
\end{corollary}

In the case of linear permutations, explicit expressions for  $P_{\ul{123}}(u,z)$, $P_{\ul{123}}(0,z)$ and $P_{\ul{1234}}(0,z)$ have been given in \cite[Theorems 4.1 and 4.3]{EliNoy1}. 
Let us also point out that, for $\si=\ul{12\dots m}$, the generating function 
$D_{\sigma}=C_{\sigma}'$ coincides with the generating function denoted by $R$ in the proof of \cite[Theorem 3.1]{EliNoy1}.

\subsection{Non-overlapping patterns}

A consecutive pattern $\sigma$ of length $m$ is called {\em non-overlapping} if two occurrences of $\sigma$ cannot overlap in more than one position; in other words, there is no permutation $\pi\in\S_{2m-2}$ with $o_\sigma(\pi)\ge 2$.

Generalizing  \cite[Theorem 3.2]{EliNoy1}, it is shown in  \cite[Theorem 3.1]{EliNoy2} that, for any non-overlapping consecutive pattern $\si$ of length $m\ge3$ with $\sigma_1=1$, the function $\omega_\si=\omega_\si(u,z)$ satisfies the following differential equation, where $b=\sigma_m$:
\begin{equation}\label{eq:nonover}\omega_\si^{(b)}+(1-u)\frac{z^{m-b}}{(m-b)!}\,\omega_\si'=0,\end{equation}
with initial conditions $\omega_\si(u,0)=1$, $\omega_\si'(u,0)=-1$, and $\omega_\si^{(i)}(u,0)=0$ for $2\le i\le b-1$.
Again, by Theorem~\ref{thm:main}, this determines $C_\si=1-\ln \omega_\si$ for all such patterns.
In this case, the generating function 
$C_{\sigma}'$ coincides with the generating function denoted by $R$ in the proof of \cite[Theorem 3.2]{EliNoy1}.

In the case $b=2$, rewriting~\eqref{eq:nonover} as a differential equation for $C_\sigma$ using~\eqref{eq:omC} and its derivatives, we obtain the following.

\begin{corollary}\label{cor:nonover}
Let $\sigma$ be a non-overlapping pattern of length $m\ge3$ with $\sigma_1=1$ and $\sigma_m=2$, and let
$D=D_{\si}(u,z)=C_{\si}'(u,z)$. Then $D$ satisfies the differential equation
\begin{equation}\label{eq:Dnonover} 
D'=D^2+(u-1)\frac{z^{m-2}}{(m-2)!}\,D
\end{equation}
with initial condition $D(u,0)=1$.
An explicit expression is given by
$$D_{\si}(u,z)=\frac{e^{(u-1)\frac{z^{m-1}}{(m-1)!}}}{1-\int_0^z e^{(u-1)\frac{t^{m-1}}{(m-1)!}}\,dt},$$
or equivalently,
$$C_{\si}(u,z)=1-\ln\left(1-\int_0^z e^{(u-1)\frac{t^{m-1}}{(m-1)!}}\,dt\right).$$
\end{corollary}

Setting $u=0$ in Equation~\eqref{eq:Dnonover} for $m=3$ gives the equation $$D_{\ul{132}}'(0,z)=D_{\ul{132}}(0,z)^2-zD_{\ul{132}}(0,z).$$ Dividing both sides by $D_{\ul{132}}(0,z)$, integrating, and using that $D_{\ul{132}}=C'_{\ul{132}}$, we obtain 
$\ln C'_{\ul{132}}(0,z)=C_{\ul{132}}(0,z)-z^2/2$, or equivalently, 
\begin{equation}\label{eq:conj2} 
C'_{\ul{132}}(0,z)=e^{C_{\ul{132}}(0,z)-z^2/2},
\end{equation}
proving part 2 of \cite[Conjecture 6.4]{DLMSSS}. 

\subsection{Other patterns and future work}

In~\cite{EliNoy2}, differential equations are also given for $\omega_\sigma(u,z)$ when $\sigma$ is any of $\ul{1324}$, $\ul{12534}$, or $\ul{13254}$. For each of these patterns, Theorem~\ref{thm:main} can again be applied to obtain $C_\sigma(u,z)$.

A natural problem for further research would be to find $C_\sigma(u,z)$ for consecutive patterns $\sigma$ that do not begin with $1$ (even after applying the basic symmetries).

In a different direction, it is shown in~\cite{EliCMP} that, for $n$ large enough, the number of (linear) permutations in $\S_n$ that avoid a consecutive pattern $\sigma$ of length $m$ is largest when $\sigma$ is a monotone pattern, and it is smallest when $\sigma=\ul{12\dots (m-2)m(m-1)}$ (or any of its symmetries). One could ask if there is an analogue of this theorem for consecutive patterns in circular permutations.


\begin{thebibliography}{}

\bibitem{ArcEli} Kassie Archer and Sergi Elizalde, Cyclic permutations realized by signed shifts, {\it J. Comb.} 5 (2014), 1--30.

\bibitem{BS} Eric Babson and Einar Steingr\'{\i}msson, Generalized permutation patterns and a classification
of the Mahonian statistics, {\it S\'em. Lothar. Combin.} 44 (2000), Art. B44b, 18 pp. 

\bibitem{Callan} David Callan, Pattern avoidance in circular permutations, preprint, \texttt{arXiv:0210014}.

\bibitem{DLMSSS} Rachel Domagalski, Jinting Liang, Quinn Minnich, Bruce E. Sagan, Jamie Schmidt and Alexander Sietsema,
Cyclic Pattern Containment and Avoidance, preprint, \texttt{arXiv:2106.02534}.

\bibitem{DK} Vladimir Dotsenko and Anton Khoroshkin, Shuffle algebras, homology, and consecutive pattern avoidance,
Algebra Number Theory 7 (2013), 673–700.

\bibitem{EliCMP} Sergi Elizalde, The most and the least avoided consecutive patterns, {\it Proc. Lond. Math. Soc.} 106 (2013), 957--979.

\bibitem{Elisurvey} Sergi Elizalde,  A survey of consecutive patterns in permutations,  Chapter in {\it Recent Trends in Combinatorics (IMA Volume in Mathematics and its Applications)}, Springer, 2016.

\bibitem{Elicont}  Sergi Elizalde, Continued fractions for permutation statistics, {\it Discrete Math. Theor. Comput. Sci.} 19 (2017), \#11.


\bibitem{EliNoy1} Sergi Elizalde and Marc Noy, Consecutive patterns in permutations, {\it Adv. in Appl. Math.} 30 (2003), 110--125.

\bibitem{EliNoy2} Sergi Elizalde and Marc Noy, Clusters, generating functions and asymptotics for consecutive patterns in permutations, {\it Adv. in Appl. Math.} 49 (2012), 351--374.

\bibitem{EliRoi1} Sergi Elizalde and Yuval Roichman, Schur-positive sets of permutations via products of grid classes, {\it J. Algebraic Combin.} 45 (2017), 363--405.

\bibitem{EliRoi2} Sergi Elizalde and Yuval Roichman, On rotated Schur-positive sets, {\it J. Combin. Theory Ser. A} 152 (2017), 121--137. 

\bibitem{FS} Philippe Flajolet and Robert Sedgewick, {\it Analytic combinatorics}, Cambridge University Press, Cambridge, 2009.

\bibitem{Kit} Sergey Kitaev, Multi-avoidance of generalised patterns, {\it Discrete Math.} 260 (2003) 89--100.

\bibitem{LR} Jeffrey Liese and Jeffrey Remmel, Generating functions for permutations avoiding a consecutive pattern, Ann. Comb. 14 (2010), 123--141.

\bibitem{MR} Anthony Mendes and Jeffrey Remmel, Permutations and words counted by consecutive patterns, {\it Adv. in Appl. Math.} 37 (2006), 443--480.

\bibitem{sagan}
Bruce~E. Sagan,
{\em Combinatorics: {T}he art of counting}, volume 210 of {\em
  Graduate Studies in Mathematics},
American Mathematical Society, Providence, RI, 2020.

\bibitem{SS} Rodica Simion and Frank W.\ Schmidt, Restricted Permutations, {\it European J. Combin.} 6 (1985), 383--406.

\bibitem{Stanley}
Richard P.\ Stanley, {\em Enumerative Combinatorics, Vol.\ 2}, Cambridge University Press, Cambridge, 1999.

\bibitem{Vella} Antoine Vella, Pattern avoidance in permutations: linear and cyclic
orders, {\it Electron. J. Combin.} 9 (2002-3), \#R18.

\end{thebibliography}
\end{document}